\newcommand{\ncm}{\newcommand}
\ncm{\aut}{auto\-mor\-phi\-sm}
\ncm{\Inn}{\mbox{\rm Inn($A$)}}
\ncm{\Ap}{\mbox{$\overline{\rm Inn}(A)$}}
\ncm{\Ext}{\mbox{\rm Ext}}
\ncm{\OExt}{\mbox{\rm OrderExt}} 
\ncm{\AI}{\mbox{\rm AInn($A$)}}
\ncm{\HI}{\mbox{\rm HInn($A$)}}
\ncm{\Aut}{\mbox{\rm Aut($A$)}}
\ncm{\Mal}{\mbox{$M_{\alpha}$}}
\ncm{\Aff}{\mbox{${\rm Aff}(T_A)$}}
\ncm{\id}{\mbox{\rm id}}
\ncm{\BE}{\begin{eqnarray*}}
\ncm{\EE}{\end{eqnarray*}}
\ncm{\lra}{\mbox{$\longrightarrow$}}
\ncm{\Hom}{\mbox{\rm Hom}}
\ncm{\el}{\ell}

\ncm{\cstar}{$C^{*}$-algebra}
\ncm{\cstars}{$C^{*}$-algebras}
\ncm{\ra}{\mbox{$\rightarrow$}}
\ncm{\al}{\mbox{$\alpha $}}
\ncm{\del}{\mbox{$\delta$}}
\ncm{\supp}{\mbox{\rm supp}}
\ncm{\Ad}{\mbox{\rm Ad}}
\ncm{\CAR}{\mbox{$M_{2^{\infty}}$}}
\ncm{\ep}{\mbox{$\epsilon > 0$}}
\ncm{\mod}{\mbox{\rm mod}}
\ncm{\ol}{\overline}
\ncm{\Mninf}{\mbox{$M_{n^{\infty}}$}}
\ncm{\MR}{M. R\o{}rdam}
\ncm{\Range}{\mbox{\rm Range}}
\ncm{\vo}{\bf}
\ncm{\ch}{\it}
\ncm{\CMP}{Comm. Math. Phys.}
\ncm{\add}{}

\newtheorem{theo}{Theorem}

\newtheorem{lem}[theo]{Lemma}
\newtheorem{prop}[theo]{Proposition}
\newtheorem{remark}[theo]{Remark}
\newtheorem{definition}[theo]{Definition}
\newtheorem{example}[theo]{Example}

\newenvironment{rem}{\begin{remark} \rm}{\end{remark}}
\newenvironment{pf}{{\it Proof.}}{\vspace{3mm}}

\newenvironment{ex}{\begin{example} \rm}{\end{example}}

\ncm{\R}{\mbox{\bf R}}
\ncm{\Z}{\mbox{\bf Z}}
\ncm{\T}{\mbox{\bf T}}
\ncm{\TT}{\T$^{2}$}
\ncm{\N}{\mbox{\bf N}}
\ncm{\C}{\mbox{\bf C}}
\documentstyle{article} 
\title{The Ext class of an approximately inner automorphism, II}
\author{A. Kishimoto and A. Kumjian}
\date{May 1998}
\begin{document}
\maketitle

\begin{abstract}
Let $A$ be a simple unital AT algebra of real rank zero and Inn($A$)
the group of inner automorphisms of $A$. In the previous paper we have
shown that the natural map of the group \Ap\ of approximately inner
automorphisms into 
$$\Ext(K_1(A),K_0(A))\oplus\Ext(K_0(A),K_1(A))
$$ 
is surjective; the kernel of this map includes the subgroup of automorphisms
which are homotopic to \Inn. In this paper we consider the quotient
of \Ap\ by the smaller normal subgroup \AI\ which consists of asymptotically
inner automorphisms and describe it as 
$$\OExt(K_1(A),K_0(A))\oplus \Ext(K_0(A),K_1(A)),
$$
where $\OExt(K_1(A),K_0(A))$ is a kind of extension group 
which takes into account the fact that $K_0(A)$ is an
ordered group and has the usual Ext as a quotient.
\end{abstract}

\section{Introduction}
An \aut\ \al\ of a unital \cstar\ $A$ is called inner if there is a
unitary $u\in A$ such that
$\al(a)=\Ad u(a)=uau^*,\ a\in A$. We denote by \Inn\ the group of
inner \aut s of $A$, which is a normal subgroup of the group
\Aut\ of all \aut s of $A$. The topology on \Aut\ is determined by
the pointwise convergence on $A$. The closure  \Ap\ of \Inn\ in
\Aut\ is, by definition, the group of approximately inner \aut s.

There are two distinguished normal subgroups of \Ap\ containing \Inn.
One is the group \HI\ of \aut s which are homotopic to \Inn,
i.e., $\al\in\HI$ if and only if there is a continuous map 
$ \al. :[0,1]\ra\Ap$ such that
$$\al_{0}\in \Inn, \ \al_1=\al.
$$
The other is the group \AI\ of asymptotically inner \aut s, i.e.,
$\al\in\AI$ if and only if there is a continuous map
$\al. :[0,1]\ra \Ap$ and a continuous map
$u. : [0,1)\ra U(A)$ with $U(A)$ the unitary group of $A$
such that
$$\al_t=\Ad u_t\ {\rm for} \ t\in [0,1),\ \al_1=\al.
$$
It is easy to show that they are indeed normal subgroups and that
$$ \Inn \subset \AI \subset \HI \subset \Ap.
$$
In this paper we describe the quotient
$$
 \Ap/\AI
$$
in terms of K-theoretic data when $A$ is a simple unital AT algebra of
real rank zero.

Recall that a unital \cstar\ $A$ is said to be a unital AT algebra if it
is expressible as the inductive limit of T algebras, i.e., finite direct
sums of matrix algebras over $C(\T)$, with unital embeddings. Note that
a unital AT algebra $A$ is stably finite and we denote by $T_A$ the convex
set of tracial states of $A$.

Let $A$ be a simple unital AT algebra of real rank zero and $\al\in \Ap$.
(In this case $\al\in \Aut$ belongs to \Ap\ if and only if $\al_*=\id$
on $K_*(A)$ \cite{El}.) 
The mapping torus of \al\ is the \cstar:
$$ M_{\alpha}=\{ x\in C[0,1]\otimes A;\ \al(x(0))=x(1)\,\}.
$$
The suspension of $A$, $SA$, is identified with the ideal of $M_{\alpha}$:
$$
SA=\{ x\in C[0,1]\otimes A;\ x(0)=0=x(1)\, \}.
$$
>From the short exact sequence:
$$
0\lra SA \lra M_{\alpha} \lra A\lra 0,
$$
one obtains the usual six-term exact sequence in K theory, which, since
$\al\in \Ap$, splits into two short exact sequences:
$$
0\lra K_i(A) \lra K_{i+1}(M_{\alpha})\lra K_{i+1}(A)\lra 0
$$
for $i=0,1$, where $K_{i+1}(SA)$ has been identified with $K_i(A)$.
Let $\eta_i(\al)$ denote the class of this sequence in
$\Ext(K_{i+1},K_i(A))$ and let $\eta$ denote the map of \Ap\ into
$$
\oplus_{i=0}^{1}\Ext(K_{i+1}(A), K_i(A))
$$
defined by $\al\mapsto (\eta_0(\al),\eta_1(\al))$, which is a group
homomorphism. (By using KK theory and the universal coefficient theorem
\cite{RS}, $\eta(\al)$ is also described as $KK(\al)-KK(\id)$.) 
In the previous paper \cite{KK1} we showed that $\eta$ induces
a surjective homomorphism:
$$
\Ap/\HI \longrightarrow \Ext(K_1(A),K_0(A))\oplus\Ext(K_0(A),K_1(A)).
$$

To state our main result of this paper we proceed to describe a
natural map $R_{\alpha}$ of $K_1(\Mal)$ into \Aff, which is the real
Banach space of affine continuous functions on the compact tracial
state space $T_A$ of $A$. Note that, since we assume that $A$ has 
real rank zero, $T_A$ is isomorphic to the state space
of $(K_0(A), [1])$. If $u\in \Mal$ is a unitary given as a piecewise
smooth function
of $[0,1]$ into $A$, then $R_{\alpha}([u])$ is defined by
$$
R_{\alpha}([u])(\tau)=\frac{1}{2\pi i} \int_0^1\tau(\dot{u}(t)u(t)^*)dt
$$
for $\tau\in T_A$. The map $R_{\alpha}$ is a group homomorphism
of $K_1(A)$ into \Aff\ and extends the natural map $D$ of $K_0(A)$
into \Aff\ when $K_0(A)$ is regarded as a subgroup of $K_1(\Mal)$.

We take the set of pairs $(E,R)$ where $E$ is an abelian group such that
$$
0 \lra K_0(A) \stackrel{\iota}{\lra} E \stackrel{q}{\lra} K_1(A) \lra 0
$$
and $R$ is a homomorphism:
$$
R: E \lra \Aff
$$
such that $R\circ\iota=D$. We can form a group $\OExt(K_1(A),K_0(A))$
from this set in much the same way as we do $\Ext(K_1(A),K_0(A))$
from the set of $E$ alone. From the previous paragraph we can associate
$\tilde{\eta}_0(\al)\in\OExt(K_1(A),K_0(A))$ with each
$\al\in\Ap$ and show that $\tilde{\eta}_0$ is a homomorphism.
Our main result is
$$\Ap/\AI \cong \OExt(K_1(A),K_0(A))\oplus \Ext(K_0(A),K_1(A))
$$
where the isomorphism is induced by the map
$\al\mapsto (\tilde{\eta}_0(\al),\eta_1(\al))$
(see Theorem~\ref{4.4}).

In Section 2 we will define $\OExt(K_1(A),K_0(A))$ and the homomorphism
$$\tilde{\eta}: \Ap\ra \OExt(K_1(A),K_0(A))\oplus \Ext(K_0(A),K_1(A))$$
in details and in Section 3 we will show that
$$
\ker\tilde{\eta}=\AI.
$$
In Section 4 we will show that $\tilde{\eta}$ is surjective; thus proving
the main result.

The authors are indebted to G.A. Elliott for discussions at an early stage
of this work. 

\section{OrderExt}
\setcounter{theo}{0}
Let $A$ be a simple unital \cstar\ and let $T_A$ be the set of tracial
states of $A$. Let $\al\in\Ap$ and let \Mal\ be the mapping torus of
\al. For a unitary $u\in\Mal$ such that $t\mapsto u(t)$ is (piecewise)
$C^1$
and for $\tau\in T_A$, we define
$$
\tau(u)=\frac{1}{2\pi i}\int_0^1\tau(\dot{u}(t)u(t)^*)dt.
$$
Since $\tau(\dot{u}(t)u(t)^*)=-\tau(u(t)\dot{u}(t)^*)$, it follows that
$\tau(u)\in\R$. If $u,v\in\Mal$ are $C^1$-unitaries, we obtain that
$$
\tau(uv)=\tau(u)+\tau(v).
$$
If $h=h^*\in \Mal$ is $C^1$, then we have for $u=e^{2\pi ih}$
$$
\tau(u)=\int_0^1\tau(\dot{h}(t))dt=\tau(h(1))-\tau(h(0))=0,
$$
where we have used that $\tau\circ\al=\tau$, which follows since
$\al\in\Ap$. Thus it follows that $\tau(u)$ is constant on each
connected component of the $C^1$-unitary group of \Mal.
By taking the matrix algebras over \Mal\ and using the density
of $C^1$-unitaries in the unitary group, we obtain a
homomorphism $\tau: K_1(\Mal)\ra \R$ by $[u]\mapsto \tau(u)$
for each $\tau\in T_A$. Since $\tau\in T_A\ra \tau(u)$ is affine
and continuous, we thus obtain:

\begin{lem}\label{2.1}
For $\al\in\Ap$ there exists a homomorphism
$$
R_{\alpha}:K_1(\Mal)\lra \Aff
$$
by $R_{\alpha}([u])(\tau)=\tau(u)$, which will be called
the {\em rotation map} for \al.
\end{lem}

Since $\al_*=\id$ on $K_i(A)$, we have the short exact sequence:
$$
0\lra K_0(A) \stackrel{\iota_*}{\lra}K_1(\Mal)\stackrel{q_*}{\lra}K_1(A)\lra 0
$$
from the short exact sequence of \cstar s:
$$
0\lra SA\stackrel{\iota}{\lra}\Mal\stackrel{q}{\lra}A\lra 0.
$$
If $p$ is a projection in $A$, we have that $\iota_*([p])=[u]$
where $u\in\Mal$ is the unitary defined by
$$
u(t)=e^{2\pi it}p +1-p.
$$
Thus we obtain:

\begin{lem}\label{2.2}
For $\al\in\Ap$ the following diagram commutes:
$$
\begin{array}{rcccl}
K_0(A)&& \stackrel{\iota_*}{\lra}&& K_1(\Mal) \\
&{\scriptstyle  D}\searrow&&\swarrow {\scriptstyle R_{\alpha}}\\
&&\Aff &&
\end{array}
$$
where $D$ is the homomorphism of $K_0(A)$ into \Aff\ defined by
$D([p])(\tau)=\tau(p)$, which will be called the {\em dimension map}
for $A$.
\end{lem}

Let $G_i=K_i(A)$. If
$$
0\lra G_0\stackrel{\iota}{\lra}E\stackrel{q}{\lra}G_1\lra0
$$
is exact, we denote this short exact sequence by $E$, the same symbol at the 
middle. Let $R$ be a homomorphism of $E$ into \Aff\ such that
$R\circ\iota=D$. We consider the set of all pairs $(E,R)$, which we
call orderextensions for $(G_1,G_0)$.

If $(E',R')$ is another orderextension, we say that $(E,R)$ and
$(E',R')$ are isomorphic with each other if there is an isomorphism
$\varphi$ of $E$ into $E'$ such that $R=R'\circ\varphi$ and
$$
\begin{array}{ccccccccc}
0 &\lra &G_0 &\stackrel{\iota}{\lra} &E &\stackrel{q}{\lra} &G_1 &\lra &0 \\
& &\parallel& &\downarrow\varphi&&\parallel&&\\
0 &\lra &G_0 &\stackrel{\iota'}{\lra}&E'&\stackrel{q'}{\lra}&G_1&\lra&0
\end{array}
$$
is commutative. Note that if $(E,R)$ and $(E',R')$ are isomorphic,
$E$ and $E'$ are isomorphic as extensions. We define an addition
for such pairs by extending that for extensions as follows. 
If $(E,R)$ and $(E',R')$ are given,
define
\begin{eqnarray*}
E''&=& \{(x,y)\in E\oplus E'| q(x)=q'(y)\}/\{(\iota(a),-\iota'(a))| a\in G_0\}\\
\iota''&:& G_0\lra E''\\
&& a\longmapsto [(\iota(a),0)]\\
q''&:& E''\lra G_1\\
&& [(x,y)]\mapsto q(x)\\
R''&:& E''\lra \Aff \\
&& [(x,y)]\mapsto R(x)+R'(y).
\end{eqnarray*}
It is easy to show that these objects are well-defined,
$$
0\lra G_0\stackrel{\iota''}{\lra}E''\stackrel{q''}{\lra}G_1\lra 0
$$
is exact, and $R''\circ\iota''=D$.
The sum of $(E,R)$ and $(E',R')$ is defined to be $(E'',R'')$. Again 
it is easy to show that the isomorphism classes of those orderextensions
form an abelian semigroup. Then the identity element for
this semigroup is given by the isomorphism class $[(E_0,R_0)]$ of
the trivial orderextension $(E_0,R_0)$ given by:
\begin{eqnarray*}
E_0&=&G_0\oplus G_1\\
\iota_0&:& G_0\lra E_0 \\
&& a\longmapsto (a,0)\\
q_0&:& E_0\lra G_1\\
&& (a,b)\mapsto b\\
R_0&:& E_0\lra \Aff\\
&&(a,b)\mapsto D(a).
\end{eqnarray*}
The inverse of $[(E,R)]$ is given by $[(E',R')]$ where
\begin{eqnarray*}
E'&=&E\\
\iota'&=&-\iota\\
q'&=& q\\
R'&=&-R.
\end{eqnarray*}
Thus the semigroup is a group, which we denote by 
$$\OExt(G_1,G_0).
$$ Note that
$\OExt(G_1,G_0)$ depends also on the dimension map $D:G_0\ra \Aff$.

\begin{lem}\label{2.3} 
The map 
\begin{eqnarray*}
\tilde{\eta}_0 &:& \Ap \lra \OExt(K_1(A),K_0(A))\\
&& \al \longmapsto [(K_1(\Mal),R_{\alpha})]
\end{eqnarray*}
is a homomorphism.
\end{lem}
\begin{pf}
By Lemma~\ref{2.2} $\tilde{\eta}_0$ is well-defined.

Let $\al,\beta\in\Ap$ and $(E,R)$ be the sum of $(K_1(\Mal),R_{\alpha})$
and $(K_1(M_{\beta}),R_{\beta})$. We have to show that $(E,R)$ is
isomorphic to $(K_1(M_{\alpha\beta}),R_{\alpha\beta})$.

Let $g\in K_1(\Mal)$ and $h\in K_1(M_{\beta})$ such that
$q(g)=q(h)$. Let $v\in M_n\otimes \Mal$ and
$w\in M_n\otimes M_{\beta}$ be unitaries such that
$[v]=g,\ [w]=h$, and $v(0)=w(0)$. Then we define a unitary
$u\in M_n\otimes M_{\alpha\beta}$ by
$$
u(t)=\left\{ \begin{array}{ll}
                 v(2t) & 0\leq t\leq 1/2\\
                \al(w(2t-1))& 1/2\leq t\leq 1
            \end{array}
\right.
$$
Then $[u]\in K_1(M_{\alpha\beta})$ depends only on $[v]$ and $[w]$.
Thus we have a map $\varphi$ of
$$
\{ (g,h)\in K_1(\Mal)\oplus K_1(M_{\beta})\,|\, q(g)=q(h)\}
$$
to $K_1(M_{\alpha\beta})$. It is easy to show that
$\varphi$ is a surjective homomorphism and the kernel of $\varphi$
equals $\{(\iota(a),-\iota(a))\,|\, a\in K_0(A)\}$. Hence $\varphi$
induces an isomorphism $\phi:E\ra K_1(M_{\alpha\beta})$.
Since
$$
R_{\alpha}([u])=R_{\alpha}([v])+R_{\beta}([w])
$$
for the above $u$, $(E,R)$ is isomorphic to
$(K_1(M_{\alpha\beta}),R_{\alpha\beta})$.
\end{pf}

\begin{lem}\label{2.4}
If $(E,R)$ is an orderextension for $(G_1,G_0)$ and $\Range\, R=\Range\, D$,
then
$$
0\lra \ker D\stackrel{\iota_*|\ker D}{\lra}\ker R\stackrel{q_*|\ker R}{\lra}
G_1\lra 0
$$
is exact.
\end{lem}
\begin{pf} 
It is obvious that the above sequence is well-defined, the compositions
of two consecutive maps vanish, and it is exact at $\ker D$. Let $g\in\ker R$
with $q_*(g)=0$. Then there is a $g'\in G_0$ such that $\iota_*(g')=g$.
But, since $D(g')=R(g)=0$, we have that $g'\in \ker D$, which implies
that it is exact at $\ker R$. Let $g\in G_1$. Then there is a $g'\in E$
with $q_*(g')=g$ and there must be a $g''\in G_0$ such that $D(g'')=R(g')$.
Since $q_*(g'-\iota_*(g''))=g$ and $R(g'-\iota_*(g''))=0$, we have that
$g\in\Range (q_*|\ker R)$.
\end{pf}

\begin{prop}\label{2.5}
If $(E,R)$ is an orderextension for $(G_1,G_0)$, the following conditions
are equivalent:
\begin{enumerate}
\item $[(E,R)]=0$,
\item \begin{enumerate}
       \item $0\ra G_0\ra E\ra G_1\ra 0$ is trivial,
        \item $\Range\, R=\Range\, D$,
        \item $0\ra \ker D\ra \ker R \ra G_1\ra 0$ is trivial,
       \end{enumerate}
\item $0\ra \ker D\ra \ker R\ra G_1\ra 0$ is exact and trivial.
\end{enumerate}
\end{prop}
\begin{pf}
If $(E_0,R_0)$ is the trivial orderextension, it satisfies (2). Any 
orderextension isomorphic to $(E_0,R_0)$ also satisfies (2). Thus (1)
implies (2).

Suppose that $(E,R)$ satisfies (2). Note that the sequence in (c)
is exact by \ref{2.4}. By (c) there is a homomorphism $\nu$ of
$G_1$ into $\ker R$ such that $q\circ\nu=\id$. Hence 
$E=\iota(G_0)\oplus \nu(G_1)$ and $R$ is given by
\begin{eqnarray*}
\iota(G_0)\oplus \nu(G_1)&\ra&\Aff\\
a+b &\mapsto& D(a).
\end{eqnarray*}
Thus $(E,R)$ is isomorphic to the trivial orderextension, i.e., (2)
implies (1).

It follows from \ref{2.4} that (2) inplies (3). The converse also
follows from the arguments in the previous paragraph.
\end{pf}

\begin{rem}\label{2.6}
By the Thom isomorphism \cite{Connes}, $K_i(\Mal)$ is isomorphic to
$K_{i+1}(A\times_{\alpha}\Z)$ as abelian groups. By extending $\tau\in T_A$
to a tracial state of $A\times_{\alpha}\Z$ and defining a natural map
$D_{\alpha}: K_0(A\times_{\alpha}\Z)\ra \Aff$, it follows that
$(K_1(\Mal),R_{\alpha})$ is isomorphic to 
$(K_0(A\times_{\alpha}\Z),D_{\alpha})$ \cite{Connes}.
See also \cite{dlHS,Pim,Bl}.
\end{rem}

\section{Asymptotically inner automorphisms}
\setcounter{theo}{0}
From now on we will assume that the \cstar\ $A$ is a simple unital
AT algebra of real rank zero. In this case by Elliott's result \cite{El}
$A$ is determined by $(K_0(A),[1],K_1(A))$ up to isomorphism,
where $K_0(A)$ is a dimension group, $K_1(A)$ is a torsion-free
abelian group, and $[1]\in K_0(A)^+$. Note that the tracial state
space $T_A$ of $A$ is identified with
the compact convex set of order-preserving homomorphisms $f: K_0(A)\ra \R$
with $f([1])=1$.

Let $\al\in\Ap$. We recall that \al\ is {\em asymptotically inner} if there
exists a continuous map $v: [0,1)\ra U(A)$ such that
$$
\al(a)=\lim_{t\ra1}\Ad\, v_t(a),\ a\in A.
$$
We denote by \AI\ the group of asymptotically inner automorphisms of $A$.
We also recall that $\tilde{\eta}$ is the homomorphism of \Ap\
into 
$$\OExt(K_1(A),K_0(A))\oplus \Ext(K_0(A),K_1(A))
$$
defined by $\al\mapsto \tilde{\eta}_0(\al)\oplus\eta_1(\al)$.

Before stating the main theorem of this section, let us recall the notion
of Bott element for pairs of almost commuting unitaries in a unital
\cstar\ $A$ \cite{EL,KK1}:  Given $u,v\in U(A)$ with $[u,v]\equiv uv-vu
\approx0$,
we associate $B(u,v)\in K_0(A)$, which is the equivalence class
of a projection close to the image of the Bott projection in 
$M_2\otimes C(\T^2)$
under the {\em quasi-homomorphism} from $M_2\otimes C(\T^2)$ into
$M_2\otimes A$ mapping the two canonical unitaries of $C(\T^2)$ into $u,v$
respectively. If $A=M_n$, this can also be given by
$$
B(u,v)=\frac{1}{2\pi i}{\rm Tr}(\log vuv^*u^*)\in \Z=K_0(M_n),
$$
where log is the logarithm
with values in $\{z;{\rm Im}(z)\in (-\pi,\pi)\}$. 
(That $B(u,v)$ is an integer follows from
the fact that the determinant of $vuv^*u^*$ is 1.) We note that 
$B(u,v)$ is invariant under homotopy of pairs of almost commuting unitaries
and that $B(u,v)=-B(u^*,v)=-B(v,u),\ B(u,v_1v_2)=B(u,v_1)+B(u,v_2)$.
We quote \cite{BEEK} for another characterization of the Bott element, which is
used to prove the following result we will need later: If
$A$ is a simple unital AT algebra of real rank zero and $u,v\in U(A)$
satisfy that $[u,v]\approx 0,\ B(u,v)=0$, 
${\rm Sp}(v)$ is almost dense in \T, and $[u]=0$, then there is
a path $u_t,\ t\in [0,1]$ in $U(A)$ such that
$[u_t,v]\approx0,\ u_0=1$, and $u_1=u$.
 
\begin{theo}\label{3.1}
Let $A$ be a simple unital AT algebra of real rank zero and let $\al\in\Ap$.
Then the following conditions are equivalent:
\begin{enumerate}
\item $\tilde{\eta}(\al)=0$,
\item $\al\in\AI$.
\end{enumerate}
\end{theo}
{\em Proof of (2) $\Rightarrow$ (1)}

Since $\eta$ is homotopy invariant, $\eta(\al)=(\eta_0(\al),\eta_1(\al))=0$ in
$\Ext(K_1(A),K_0(A))\oplus \Ext(K_0(A),K_1(A))$.

We may suppose that we have a piecewise $C^1$ map $v$ of $[0,1)$ into $U(A)$
such that
$$
\al(a)=\lim_{t\ra1}\Ad\,v_t(a),\ a\in A.
$$
Let $u\in U(A)$. We define a unitary $\hat{u}\in\Mal\otimes M_2$
by composing the following paths:
$$
[0,1]\ni t\mapsto R_t\left(\begin{array}{cc}1&0\\ 0& v_0\end{array}\right)
   R_t^{-1}\left(\begin{array}{cc}u&0\\0&1\end{array}\right) R_t
   \left(\begin{array}{cc}1&0\\0&v_0^*\end{array}\right)R_t^{-1}
$$
and
$$
[0,1)\ni t\mapsto \left(\begin{array}{cc}v_tuv_t^*&0\\0&1\end{array}\right)
$$
with
$$
1\mapsto \left(\begin{array}{cc}\al(u)&0\\0&1\end{array}\right),
$$
where
$$
R_t=\left(\begin{array}{cc}\cos \frac{\pi}{2}t& -\sin\frac{\pi}{2}t\\
     \sin\frac{\pi}{2}t&\cos\frac{\pi}{2}t\end{array}\right).
$$
Then it follows that
$\tau(\dot{\hat{u}}(t)\hat{u}(t)^*)=0$ for $\tau\in T_A$. In particular
$R_{\alpha}([\hat u])=0$.
Since $q_*([\hat u])=[u]$, the map $[u]\mapsto [\hat u]$ defines
a homomorphism $\varphi$ of $K_1(A)$ into $\ker R_{\alpha}$ such that
$q_*\circ\varphi=\id$. This implies that
$$
0\lra \ker D\lra \ker R_{\alpha}\lra K_1(A)\lra 0
$$
is exact and trivial, and thus concludes the proof by \ref{2.5}.
\medskip

The rest of this section will be devoted to the proof of (1)$\Rightarrow$(2).

Let $\{A_n\}$ be an increasing sequence of T subalgebras of $A$ such that
$A_1\ni 1$ and $A=\ol{\cup_{n=1}^{\infty}A_n}$. We express $A_n$ as
$$
A_n=\oplus_{i=1}^{k_n}B_{n,i}\otimes C(\T)
$$
where $B_{n,i}$ is isomorphic to the full matrix algebra $M_{[n,i]}$.
By identifying $K_i(A)$ with $\Z^{k_n}$ in a natural way we obtain a
homomorphism $K_i(A_n)$ into $K_i(A_{n+1})$ as the multiplication of 
a matrix $\chi_n^i$. We always assume that $\chi_n^0(i,j)$ is big
and $|\chi_n^1(i,j)|/\chi_n^0(i,j)$ is small compared with 1 and that
the embedding of $A_n$ into $A_{n+1}$ is in a standard form, i.e.,
$B_n=\oplus_{i=1}^{k_n}B_{ni}\subset B_{n+1}$ and the canonical unitary
$z_n$ of $1\otimes C(\T)\subset A_n$ in $B_{n+1}\cap B_n'\otimes C(\T)$ is a
direct sum of elements of the form:
$$
\left(\begin{array}{cccc} 0&&&z_{n+1}^L\\ 1&\cdot&&\\ &\ddots&\ddots&\\
  &&1&0 \end{array}\right)
$$
with $L=\pm1$; e.g., if $\chi_n^1(i,j)>0$, $z_np_{n+1\,i}p_{nj}$ is a direct
sum of $\chi_n^1(i,j)$ matrices of the above form with $L=1$ in
$B_{n+1}\cap B_n'\otimes C(\T)p_{n+1\, i}p_{nj}
\cong M_{\chi_n^0(i,j)}\otimes C(\T)$ \cite{El,KK1}. 

For each $n=1,2,\ldots$ let
$$
M_{\alpha,n}=\{x\in C[0,1]\otimes A\, |\,x(0)\in A_n,\ \al(x(0))=x(1)\}.
$$
Then we obtain the exact sequence of \cstar s:
$$
0\lra SA \stackrel{\iota_n}{\lra}M_{\alpha,n}\stackrel{q_n}{\lra}A_n\lra0
$$
from which follow the exact sequences of abelian groups:
$$
0\lra K_i(A)\lra K_{i+1}(M_{\alpha,n})\lra K_{i+1}(A_n)\lra 0.
$$
Since $K_i(A_n)\cong \Z^{k_n}$, the above extensions are all trivial.

Let $R=R_{\alpha}$ and $R_n=R\circ j_{n*}: K_1(M_{\alpha,n})\ra \Aff$,
where $j_n$ is the embedding of $M_{\alpha,n}$ into \Mal.
Since $\Range\, D=\Range\, R_n$, we obtain by \ref{2.4} that
$$
0\lra \ker D\stackrel{\iota_{n*}}{\lra}\ker R_n\stackrel{q_{n*}}{\lra}
K_1(A_n)\lra 0
$$
is exact. Note that the inductive limit of these extensions is naturally
isomorphic to the exact sequence: 
$$0\ra\ker D\ra \ker R\ra K_1(A)\ra0.$$
We shall specify a homomorphism $\varphi_n$ of $K_1(A_n)$
into $\ker R_n$ such that
$$
q_{n*}\circ\varphi_n=\id.
$$

Since $\al\in \Ap$, we have a $u_n\in U(A)$ for each $n$ such that
\begin{eqnarray*}
\al|B_n&=&\Ad\, u_n|B_n,\\
\al(z_n)&\approx &\Ad\,u_n(z_n),
\end{eqnarray*}
where $B_n=\oplus_{i=1}^{k_n}B_{ni}\subset A_n$ and $z_n$ is the canonical 
unitary of $C(\T)\subset A_n$. Define
$$
h_{ni}=\frac{1}{2\pi i}\log\,\al(z_{ni})\Ad\,u_n(z_{ni})^*
$$
where $z_{ni}=z_np_{ni}+1-p_{ni}$ with $p_{ni}$ the identity of $B_{ni}$
and $h_{ni}=h_{ni}^*$ is defined uniquely as $\|h_{ni}\|\approx0$
since $\al(z_{ni})\Ad\,u_n(z_{ni}^*)\approx 1$.
Define $\zeta_{ni}\in U(M_{\alpha,n}\otimes M_2)$ by composing two
paths of unitaries:
$$
[0,1]\ni t\mapsto R_t(u_n^*\oplus1)R_t^{-1}(u_n\oplus1)(z_{ni}\oplus1)
         (u_n^*\oplus1)R_t(u_n\oplus1)R_t^{-1}
$$
and
$$
[0,1]\ni t\mapsto e^{2\pi ith_{ni}} \Ad\,u_n(z_{ni})\oplus 1.
$$
Then we have that
\begin{eqnarray*}
q_n(\zeta_{ni})&=& z_{ni}\oplus 1,\\
R_n([\zeta_{ni}])&=& \hat{h}_{ni},
\end{eqnarray*}
where $\hat{h}_{ni}\in \Aff$ is defined by
$$
\hat{h}_{ni}(\tau)=\tau(h_{ni}),\ \tau\in T_A.
$$
Since the above procedure applies to a unitary $z_np +1-p$ with
$p$ a minimal projection in $B_n$, it follows that 
$[\zeta_{ni}]\in K_1(M_{\alpha,n})$
is divisible by $[n,i]$. Thus one obtains a homomorphism $\varphi_n$
of $K_1(A_n)$ into $K_1(M_{\alpha,n})$ with $q_{n*}\circ\varphi=\id$
by setting
$$
\varphi_n:[z_{ni}]\longmapsto [\zeta_{n,i}].
$$

\begin{lem}\label{3.2}
$\Range\, D$ is dense in \Aff.
\end{lem}
\begin{pf}
Since $A$ is a simple unital AT algebra of real rank zero, it is approximately
divisible \cite{El2}. Thus this is 3.14(a) of \cite{BKR}.
(A unital \cstar\ is approximately divisible if it has a central sequence
$\{B_n\}$ of unital $C^*$-subalgebras with $B_n\cong M_2\oplus M_3$ 
\cite{BKR}. Since $A$ is obtained as the inductive limit of $\{A_n\}$
of T algebras with unital embeddings and the embeddings need to satisfy
only the K-theoretic conditions and the condition of real rank zero
\cite{BBEK}, thanks to Elliott's result \cite{El}, 
we can easily arrange the inductive system so that $A_{n+1}\cap A_n'\supset
M_2\oplus M_3$, which implies that $A$ is approximately divisible.)
\end{pf}

Let 
$$
\delta_n=\min_i\inf \{\tau(p_{ni});\ \tau\in T_A\},
$$
where $p_{ni}$ is the identity of $B_{ni}$.
Since $A$ is simple, $\delta_n$ is strictly positive. We choose the unitary 
$u_n\in A$ so that $\|h_{ni}\|<\delta_n$.
Since $\Range\, R_n=\Range\, D$, we have, for any $\epsilon>0$ with
$\|h_{ni}\|+\epsilon<\delta_n$, projections $p_{\pm}\in A$ such that
\begin{eqnarray*}
\frac{1}{[n,i]}\hat{h}_{ni}&=&D(p_+)-D(p_-),\\
\|D(p_{\pm})\|&<& \frac{1}{[n,i]}(\|h_{ni}\|+\epsilon),
\end{eqnarray*}
where $D$ is also regarded as a map of the projections into \Aff. 
(First we approximate $\hat{h}_{ni+}/[n,i]$ by $D(p_+)$ with 
$p_+$ a projection such that $D(p_+)-\hat{h}_{ni+}/[n,i]>0$ (or strictly
positive), where $h_{ni+}$ is the positive part of $h_{ni}$. We should note
that $\|\hat{h}_{ni+}/[n,i]\|\leq\|h_{ni}\|/[n,i]$
and find a projection $p_-$ such that
$D(p_-)=D(p_+)-\hat{h}_{ni}/[n,i]\approx\hat{h}_{ni-}/[n,i]$.)
Since $D(p_{\pm})<\delta_n/[n,i]\leq D(p_{ni})/[n,i]$, we find projections
$e_{i\pm}\in p_{ni}Ap_{ni}\cap B_{ni}'$ such that
\begin{eqnarray*}
\hat{h}_{ni}&=&D(e_{i+})-D(e_{i-}),\\
\|D(e_{i\pm})\|&<&\|h_{ni}\|+\epsilon.
\end{eqnarray*}
Thus, by making $\|h_{ni}\|$ small, we can make $\|D(e_{i\pm})\|$ arbitrarily
small. Then, by using Lemma~\ref{3.4} below, we can find a unitary 
$w_{ni}\in p_{ni}Ap_{ni}\cap B_{ni}'$ such that
\BE
w_{ni}&=&w_{ni}p_{ni}+1-p_{ni},\\
\Ad\, w_{ni}(z_{ni})&\approx& z_{ni},\ {\rm (in\ the\ order\ of}\ \|h_{ni}\|
    {\rm )}\\ 
\hat{k}_{ni}&=& \hat{h}_{ni},
\EE
where
$$
k_{ni}=\frac{1}{2\pi i}\log\,\Ad\,w_{ni}(z_{ni})z_{ni}^*.
$$

Let $w_n=w_{n1}w_{n2}\cdots w_{nk_n}$. Note that
\BE
\al(z_{ni})\Ad\,u_nw_n(z_{ni}^*)&=&\al(z_{ni})\Ad\,u_n(z_{ni}^*)
                \Ad\,u_n(z_{ni}\Ad\,w_n(z_{ni}^*))\\
 &=&e^{2\pi ih_{ni}}\Ad\,u_n(e^{-2\pi ik_{ni}}).
\EE
Then composing the two paths:
$$
[0,1]\ni t\longmapsto \Ad\,u_n(e^{-2\pi itk_{ni}})
$$
and
$$
[0,1]\ni t\longmapsto e^{2\pi ith_{ni}}\Ad\,u_n(e^{-2 \pi ik_{ni}})
$$
multiplied with $\Ad\,u_nw_n(z_{ni})$ to the right, we obtain a path $U$
from $\Ad\,u_nw_n(z_{ni})$ to $\al(z_{ni})$ such that
$$
\frac{1}{2\pi i}\int_0^1\tau(\dot{U}(t)U(t)^*)dt=0,\ \tau\in T_A.
$$
Since $U$ is in a small neighbourhood of $\al(z_{ni})\approx 
\Ad\,u_nw_n(z_{ni})$, it follows that the unitary $\zeta_{ni}$ obtained
from $z_{ni}$ in the same way as before with $u_nw_n$ in place of $u_n$
satisfies
$$ 
R_n([\zeta_{ni}])=0.
$$
Thus we have shown:

\begin{lem}\label{3.3}
Suppose that $\tilde{\eta}_0(\al)=0$. Then for any $n$ and $\epsilon\in(0,1)$
there exists a unitary $u_n\in A$ such that
\BE
&& \al|B_n=\Ad\,u_n|B_n,\\
&& \|\al(z_{ni})-z_{ni}\|<\epsilon,\\
&&\hat{h}_{ni}=0,
\EE
where
$$
h_{ni}=\frac{1}{2\pi i}\log\,\al(z_{ni})\Ad\,u_n(z_{ni}^*).
$$
Hence defining a unitary $\zeta_{ni}\in M_{\alpha,n}\otimes M_2$ by
composing the two paths:
$$
[0,1]\ni t\mapsto R_t(u_n^*\oplus 1)R_t^{-1}(\Ad\,u_n(z_{ni})\oplus 1)
              R_t(u_n\oplus 1)R_t^{-1}
$$
and
$$
[0,1]\ni t\mapsto e^{2\pi ith_{ni}}\Ad\,u_n(z_{ni})\oplus 1,
$$
where $R_t$ is defined as before, one can define a
homomorphism $\varphi_n$ of $K_1(A_n)$ into $\ker R_n$ by
$\varphi([z_{ni}])=[\zeta_{ni}],\ i=1,\ldots,k_n$.
\end{lem}

\begin{lem}\label{3.4}
If $e\in p_{ni}Ap_{ni}\cap B_{ni}'$ is a projection such that
$\|D(e)\|$ is sufficiently small, then
for any $\epsilon>0$ there exists a unitary
$w_{\pm}\in p_{ni}Ap_{ni}\cap B_{ni}'$ such that
\BE
&& \|\Ad\,w_{\pm}(z_{ni})-z_{ni}\|< 2\pi\|D(e)\|+\epsilon,\\
&& [w_{\pm}]=0,\\
&&B(w_{\pm},z_{ni})=\pm[e].
\EE
In particular if $k_{\pm}=\frac{1}{2\pi i}\log\,\Ad\,z_{ni}w_{\pm}(z_{ni}^*)$, 
it follows
that $\hat{k}_{\pm}=\pm D(e)$.
\end{lem}
\begin{pf}
To simplify the notation we may suppose that $p_{ni}A_{ni}p_{ni}\cap B_{ni}'$
to be $A$ and $z_{ni}$ to be the canonical unitary $z_1\in A_1=C(\T)$.

Since the projection $e$ plays a role only through $[e]$, we may suppose that
$e\in A_m$ for some $m>1$. We will later assume that $m$ is sufficiently
large. Since $A_n\hookrightarrow A_{n+1}$ are in the standard form,
$z_1p_{mj}$ in $B_{mj}\otimes C(\T)$ looks like a direct sum of elements of
the form:
$$
\left(\begin{array}{cccc}0&&&z_{mj}^{L_s}\\
       1&\cdot&&\\
       &\ddots&&\\
       &&1&0
       \end{array}\right)\in M_{M_s}(C(\T))
$$
where $L_s=\pm1,\ M_s \gg 1$ and
\BE
\sum_s L_s&=&\chi_{m1}^1(j,1),\\
\sum_s M_s&=& \chi_{m1}^0(j,1)=[m,j].
\EE
Note that $D(e)$ takes values in the convex hull of
$$
\frac{\dim(ep_{mj})}{[m,j]},\ j=1,\ldots, k_m,
$$
which are all assumed to be much less than 1. Let $t_m$ be the maximum
of these $k_m$ values. Then $t_m$ decreases as $m\ra\infty$ and the limit 
of $t_m$ equals $\tau(e)$ for some $\tau\in T_A$ (or $\|D(e)\|$).
Thus if $m$ is sufficiently large, we may assume that 
$t_m<\|D(e)\|+\epsilon/4\pi$. We can obtain the required
unitary $w_j$ in $B_{mj}\otimes C(\T)$ 
as the direct sum of elements of the form:
$$
\left(\begin{array}{ccccc}
 1&&&&\\ &\omega&&&\\&&\omega^2&&\\&&&\ddots&\\&&&&\omega^{M_s-1}
 \end{array}\right)
$$
where $\omega=e^{-2\pi iN_s/M_s}$ and the integers $N_s$ are chosen so that
\BE
\sum N_s &=& \dim(ep_{mj}),\\
\frac{N_s}{M_s}&\approx& \frac{\dim(ep_{mj})}{[m,j]}.
\EE
Note that by defining
$$
k_j=\frac{1}{2\pi i}\log\, z_1p_{mj}\Ad w_j(z_1^*p_{mj}),
$$
the Bott element $B(w_j,z_1p_{mj})\in K_0(A_mp_{mj})=\Z$ 
for the almost commuting pair $w_j,z_1p_{mj}$ of unitaries 
in $A_mp_{mj}=B_{mj}\otimes C(\T)$ is equal to
$$
 {\rm Tr}(k_j)={\rm Tr}(\oplus_s\frac{N_s}{M_s}1_s)=\sum N_s =\dim(ep_{mj}),
$$
where $k_j\in B_{mj}\otimes C(\T)$ should be evaluated at some (or any) 
point of \T\ (see \cite{EL,KK1,BEEK}).
This shows that
$$
B(w_j,z_1p_{mj})=[ep_{mj}],
$$
and in particular that $\hat{k}_j=D(ep_{mj})$.

If $m$ is sufficiently large or all $M_s$ are sufficiently large, we can 
assume that
$$
\frac{N_s}{M_s}< \|D(e)\|+\epsilon/2\pi.
$$
Thus we obtain the norm estimate
$$
\|\Ad\,w_j(z_1p_{mj})-z_1p_{mj}\|< 2\pi \|D(e)\|+\epsilon.
$$
By taking $w_+=w_1+w_2+\cdots +w_{k_m}$, this completes the proof for
$w_+$. For $w_-$ we just replace $\omega$ in the definition
of $w_j$ by $\ol{\omega}=e^{2\pi iN_s/M_s}$.
\end{pf}

By defining $\varphi_n: K_1(A_n)\ra \ker R_n$ as above, we identify
$\ker R_n$ with $\ker D\oplus K_1(A_n)$. We now have to translate
the natural map $\ker R_n\ra \ker R_{n+1}$ into the map
$\psi_n : \ker D\oplus K_1(A_n)\ra \ker D\oplus K_1(A_{n+1})$:
$$
\begin{array}{ccccccccccc}
0 &\ra& \ker D&\ra &\ker D&\oplus&K_1(A_n)&\ra&K_1(A_n)&\ra& 0\\
&& \parallel&& \parallel &\psi_n^0\swarrow&\downarrow\chi_n^1&&\downarrow\chi_n^1
&&\\
0&\ra&\ker D&\ra&\ker D&\oplus&K_1(A_{n+1})&\ra&K_1(A_{n+1})&\ra&0
\end{array}
$$
where we have used that $\psi_n$ must be of the form
$\psi_n(a,b)=(a+\psi_n^0(b),\chi_n^1(b))$.

\begin{lem}\label{3.5}
If $u_n$ is a unitary in $A$ and $\epsilon\in (0,1)$ such that
\BE
&& \al|B_n=\Ad\,u_n|B_n,\\
&& \|\al(z_n)-\Ad\,u_n(z_n)\|<\epsilon,\\
&& \hat{h}_{ni}=0,
\EE
then for any $m\leq n$ and $j=1,\ldots,k_m$,
\BE
& \|\al(z_{mj})-\Ad\,u_n(z_{mj})\|<\epsilon,&\ \ \ (*)\\
& \hat{h}_{mj}=0,&\ \ \ (**)
\EE
where
$$
h_{mj}=\frac{1}{2\pi i}\log\, \al(z_{mj})\Ad\,u_n(z_{mj}^*).
$$
\end{lem}
\begin{pf}
By the assumption on the embedding of $A_m$ into $A_n$, ($*$) follows
immediately. Since the homomorphism $\varphi_n:K_1(A_n)\ra \ker R_n$
can be defined on $[z_{mj}]$ in the canonical way and 
$R_n\varphi_n([z_{mj}])=\hat{h}_{mj}$, ($**$) also follows 
immediately.
\end{pf}

\begin{lem}\label{3.6}
The homomorphism $\psi_n^0:K_1(A_n)\ra \ker D$ is given by
$$
[z_{ni}]\longmapsto B(u_{n+1}^*u_n,z_{ni}),
$$
where $[z_{ni}]=[n,i]e_i$ with $(e_i)_i$ the canonical basis for
$\Z^{k_n}=K_1(A_n)$ and $B(u_{n+1}^*u_n,z_{ni})$ is divisible
by $[n,i]$.
\end{lem}
\begin{pf}
First of all we shall show that
$D(B(u_{n+1}^*u_n,z_{ni}))=0$. Because if we define self-adjoint
$h_i\in A$ by
\BE
 h_1&=&\frac{1}{2\pi i}\log\,\al(z_{ni})\Ad\,u_{n+1}(z_{ni}^*),\\
 h_2&=&\frac{1}{2\pi i}\log\,\al(z_{ni})\Ad\,u_n(z_{ni}^*),\\
 h_3&=&\frac{1}{2\pi i}\log\,z_{ni}\Ad(u_{n+1}^*u_n)(z_{ni}^*),
\EE
then $\hat{h}_2=0$ and $\hat{h}_1=0$ by \ref{3.5} and hence 
$\hat{h}_3=0$ since
$$
\Ad\,u_{n+1}(e^{2\pi ih_3})=e^{-2\pi ih_1}e^{2\pi ih_2}.
$$
(One way of proving that $\hat{h}_3=0$ 
is to take a closed path $w$ of unitaries: 
$$
w(t)=\left\{\begin{array}{ll}e^{-6\pi ith_1}&0\leq t\leq 1/3\\
                e^{-2\pi ih_1}e^{2\pi i(3t-1)h_2}&1/3\leq t\leq 2/3\\
             e^{-2\pi ih_1}e^{2\pi ih_2}\Ad\,u_{n+1}(e^{2\pi i(3t-2)h_3})
                & 2/3\leq t\leq 1
           \end{array}\right.
$$
in a neighbourhood of 1, and compute for any $\tau\in T_A$,
$0=1/2\pi i\int_0^1\tau(\dot{w}(t)w(t)^*)dt=-\tau(h_1)+\tau(h_2)-\tau(h_3)$.)

We may suppose that $u_{n+1}^*u_n\in A_m\cap B_n'$ for some $m>n$.
In this case $B(u_{n+1}^*u_n,z_{ni})$ in $K_0(A_m)$ is defined by
$$
\left({\rm Tr}_{B_{mj}}(h_3p_{mj})\right)_j,
$$
where $h_3p_{mj}\in B_{mj}\otimes C(\T)$ is evaluated at a point of \T\,
and $\hat h_3=0$ means that for any $\tau\in T_A$,
$$
\sum_j\tau(p_{mj})\frac{{\rm Tr}_{B_{mj}}(h_3p_{mj})}{[m,j]}=0.
$$

Define a path $v_{nt},\ t\in [0,1]$ of unitaries in $A\otimes M_2$ by
$$
v_{nt}=R_t (u_n^*\oplus 1)R_t^{-1}(u_n\oplus 1).
$$
Then to compute $\psi_n^0([z_{ni}])$ we have to calculate
$$
\psi_n^0([z_{ni}])=\varphi_n([z_{ni}])-\varphi_{n+1}([z_{ni}])=
[t\mapsto \Ad\,v_{n,t}(z_{ni})]-[t\mapsto \Ad\,v_{n+1,t}(z_{ni})]\ \ \ (*)
$$
in $K_1(M_{\alpha,n+1})$ where $z_{ni}$ is identified with
$z_{ni}\oplus 1$ (see 2.8 of \cite{KK1} for a similar computation). 
More precisely we have to add a short path
from $\Ad\,u_n(z_{ni})$ (resp. $\Ad\,u_{n+1}(z_{ni})$) to
$\al(z_{ni})$ to the path $t\mapsto \Ad\,v_{n,t}(z_{ni})$
(resp. $t\mapsto \Ad\,v_{n+1,t}(z_{ni})$) to get a unitary in
$M_{\alpha,n+1}\otimes M_2$ and we always understand the formulae
in this way. Note that ($*$) is equal to
$$
[t\mapsto \Ad\,v_{n,t}(z_{ni})\Ad\,v_{n+1,t}(z_{ni}^*)]
$$
in $K_1(SA)\subset K_1(M_{\alpha,n+1})$ or, by applying 
$t\mapsto \Ad\,v_{n+1,t}^*$, which induces the identity map on $K_1(SA)$,
to
$$
[t\mapsto v_{n+1,t}^*v_{n,t}z_{ni}v_{n,t}^*v_{n+1,t}z_{ni}^*].
$$
Since 
$$
v_{n+1,t}^*v_{n,t}=(u_{n+1}^*\oplus 1)R_t(u_{n+1}u_n^*\oplus 1)R_t^{-1}
         (u_n\oplus1),
$$
the above element is equal to the class of
$$
t\mapsto (u_{n+1}z_{ni}^*u_{n+1}^*\oplus1)R_t(u_{n+1}u_n^*\oplus1)R_t^{-1}
    (u_nz_{ni}u_n^*\oplus1)R_t(u_nu_{n+1}^*\oplus1)R_t^{-1}
$$
by applying $\Ad(u_{n+1}z_{ni}^*\oplus 1)$. Again this is equal to
the class of
$$
t\mapsto (u_n^*u_{n+1}z_{ni}^*u_{n+1}^*u_n\oplus1)R_t(u_n^*u_{n+1}\oplus1)
  R_t^{-1}(z_{ni}\oplus1)R_t(u_{n+1}^*u_n\oplus1)R_t^{-1}
$$
by applying $t\mapsto \Ad(u_n^*\oplus u_n^*)$. More precisely we have to
add a short path to connect the value at $t=1$
$$
u_n^*u_{n+1}z_{ni}^*u_{n+1}^*u_nz_{ni}\oplus 1
$$
to 1.
Since $u_{n+1}^*u_n\in A_m\cap B_n'$ by the assumption, the path can be
taken in $A_m$. The above element in $K_1(SA_m)=K_0(A_m)$ is equal to
\BE
&&\left(-\frac{1}{2\pi i}
  {\rm Tr}_{B_{mj}}\log(u_n^*u_{n+1}z_{ni}^*u_{n+1}^*u_nz_{ni}p_{mj})
\right)_j\\ 
&&=\left(\frac{1}{2\pi i}{\rm Tr_{B_{mj}}}\,\log
(z_{ni}(u_{n+1}^*u_n)z_{ni}^*(u_{n+1}^*u_n)^*p_{mj})
\right)_j\\
&&=B_{A_m}(u_{n+1}^*u_n,z_{ni}).
\EE
Note also that since the non-trivial
part of $z_{ni}(u_{n+1}^*u_n)z_{ni}^*(u_{n+1}^*u_n)^*$ belongs to
$p_{ni}A_mp_{ni}\cap B_{ni}'$, each component of 
$B_{A_m}(u_{n+1}^*u_n,z_{ni})$ is divisible by $[n,i]$. Then we obtain that
$$
\psi_n^0([z_{ni}])=B(u_{n+1}^*u_n,z_{ni}),\ \ i=1,\ldots,k_n,
$$
is a well-defined homomorphism of $K_1(A_n)$ into
$\ker D\subset K_0(A)$.
\end{pf}

\begin{lem}\label{3.7}
Suppose that $\tilde{\eta}_0(\al)=0$. Then there exist unitaries
$u_n\in A$ such that
\BE
&& \al|B_n=\Ad\,u_n|B_n,\\
&& \|\al(z_m)-\Ad\,u_n(z_m)\|<2^{-n},\ m\leq n,\\
&& B(u_{n+1}^*u_n,z_{ni})=0, \ i=1,\ldots,k_n,\\
&& \hat{h}_{ni}=0, \ \ i=1,\ldots,k_n,
\EE
where
$$
h_{ni}=\frac{1}{2\pi i}\log\,\al(z_{ni})\Ad\,u_n(z_{ni}^*).
$$
\end{lem}
\begin{pf}
By the assumption and Proposition~\ref{2.5} the sequence of trivial extensions:
$$
\begin{array}{ccccccccccc}
0&\lra&\ker D&\lra&\ker D&\oplus&K_1(A_n)&\stackrel{q_{n}}{\lra}&K_1(A_n)
   &\lra&0\\
&&\parallel&&\parallel&\psi_n^0\swarrow&\downarrow\chi_n^1&&\downarrow\chi_n^1
&&\\
0&\lra&\ker D&\lra&\ker D&\oplus&K_1(A_{n+1})&\stackrel{q_{n+1}}{\lra}&
K_1(A_{n+1})&\lra&0\\
&&\parallel&&\parallel&\swarrow&\downarrow&&\downarrow&&
\end{array}
$$
defines the trivial extension in $\Ext(K_1(A),\ker D)$. Hence
we have a homomorphism $h_n^0:K_1(A_n)\ra\ker D$ for each $n$
such that
$$
\psi_n^0=h_n^0-h_{n+1}^0\chi_n^1.
$$
(To see this we denote by $E$ the inductive limit of the middle terms,
and by $\varphi$ a homomorphism of $K_1(A)$ into $E$ such that
$q\varphi=\id$. If $\xi_n$ denotes the natural homomorphism of $K_1(A_n)$
into $\ker D\oplus K_1(A_n)$ composed with $\ker D\oplus K_1(A_n)
\ra E$, $\psi_n^0$ is given by $\psi_n^0=\xi_n-\xi_{n+1}\chi_n^1$.
We set $h_n^0=\xi_n-\varphi_n$ where $\varphi_n$ is the homomorphism
$K_1(A_n)\ra K_1(A)$ composed with $\varphi:K_1(A)\ra E$. Then it
follows that
$$
h_n^0-h_{n+1}^0\chi_n^1=\xi_n-\varphi_n-\xi_{n+1}\chi_n^1+\varphi_{n+1}\chi_n^1
=\xi_n-\xi_{n+1}\chi_n^1=\psi_n^0,
$$
where we have used that $\varphi_n=\varphi_{n+1}\chi_n^1$.)

Since $h_n^0(e_i^n)\in\ker D$, where $(e_i^n)_{i=1}^{k_n}$ is the canonical
basis for $\Z^{k_n}=K_1(A_n)$, we can find projections
$e_{i\pm}^n\in p_{ni}Ap_{ni}\cap B_{ni}'$ such that
$$
[n,i]h_n^0(e_i^n)=[e_{i+}^n]-[e_{i-}^n]
$$
and $\|D(e_{i\pm}^n)\|$ is arbitrarily small. (We find a positive 
$g\in K_0(A)$ with $\|D(g)\|$ sufficiently small and then find projections
$e_{i\pm}^n$ such that $[e_{i+}^n]=[n,i](g+h_n^0(e_i^n))$ and
$[e_{i-}^n]=[n,i]g$.) Then by Lemma~\ref{3.4} we find a unitary
$w_n\in A\cap B_n'$ such that
\BE
&& [w_n]=0,\\
&& B(w_n,z_{ni})=-[e_{i+}^n]+[e_{i-}^n]=-[n,i]h_n^0(e_i^n)
\EE
and $\|[w_n,z_{ni}]\|$ is arbitrarily small for $i=1,\ldots,k_n$. Since
\BE
B(w_{n+1}^*,z_{n,i})&=&\sum_j B(w_{n+1}^*p_{n+1,j},z_{ni}p_{n+1,j})\\
  &=&\sum_j \chi_n^1(j,i)[n,i]B(w_{n+1}^*,z_{n+1,j})/[n+1,j]\\
  &=&\sum_j\chi_n^1(j,i)[n,i]h_{n+1}^0(e_j^{n+1})\\
  &=&[n,i]h_{n+1}^0\chi_n^1(e_i^n),
\EE
we have that
$$
B(w_{n+1}^*u_{n+1}^*u_nw_n,z_{ni})=0.
$$
Since $D(B(w_n,z_{ni}))=0$, we have that $\hat{k}_i=0$ for
$k_i=1/2\pi i\log\,z_{ni}\Ad\,w_n(z_{ni}^*)$, and hence that
$\hat{h}_i=0$ for $h_i=1/2\pi i\log\,\al(z_{ni})\Ad\,u_nw_n(z_{ni}^*)$.
Thus by replacing $u_n$ by $u_nw_n$, we have the conclusion.
\end{pf}

Note that the exact sequence
$$
0\lra K_1(A)\lra K_0(\Mal)\lra K_0(A)\lra0
$$
is obtained as the inductive limit of
$$
\begin{array}{ccccccccc}
0&\lra&K_1(A)&\lra& K_0(M_{\alpha,n})&\lra&K_0(A_n)&\lra&0\\
&&\parallel&&\downarrow{\scriptstyle \psi_n}&&\downarrow
      {\scriptstyle \chi_n^0}&&\\
0&\lra&K_1(A)&\lra&K_0(M_{\alpha,n+1})&\lra&K_0(A_{n+1})&\lra&0\\
&&\parallel&&\downarrow&&\downarrow&&
\end{array}.
$$
By defining a homomorphism $\varphi_n:K_0(A_n)\ra K_0(M_{\alpha,n})$
just as in Lemma~\ref{3.3}, we identify $K_0(M_{\alpha,n})$ with
$K_1(A)\oplus K_0(A_n)$ and find a homomorphism
$\psi_n^1:K_0(A_n)\ra K_1(A)$ as in the following diagram:
$$
\begin{array}{ccccccccccc}
0&\lra&K_1(A)&\lra& K_1(A)&\oplus&K_0(A_n)&\lra&K_0(A_n) &\lra&0\\
&&\parallel&&\parallel&{\scriptstyle \psi_n^1}\swarrow
      &\downarrow{\scriptstyle \chi_n^0}&&\downarrow{\scriptstyle \chi_n^0}
&&\\
0&\lra& K_1(A)&\lra&K_1(A)&\oplus&K_0(A_{n+1})&\lra&K_0(A_{n+1})&\lra&0\\
&&\parallel&&\parallel&\swarrow&\downarrow&&\downarrow&&
\end{array}
$$
  
\begin{lem}\label{3.8}
The homomorphism $\psi_n^1:K_0(A_n)\ra K_1(A)$ is given by
$$
[p_{ni}]\mapsto [u_{n+1}^*u_np_{ni}]
$$
where $[p_{ni}]=[n,i]e_i$ with $(e_i)$ the canonical basis for
$\Z^{k_n}=K_0(A_n)$ and $[u_{n+1}^*u_np_{ni}]$ is divisible by
$[n,i]$.
\end{lem}
\begin{pf}
As in the proof of Lemma~\ref{3.6} we have to decide
$$
[t\mapsto \Ad\,v_{n,t}(p_{ni})]-[t\mapsto \Ad\,v_{n+1,t}(p_{ni})]\ \ \ (*)
$$
in $K_0(M_{\alpha,n+1})$, where $p_{ni}$ denotes $p_{ni}\oplus0$
in $A\otimes M_2$. (Note that $\Ad\,u_n(p_{ni})=\al(p_{ni})$ and
$\Ad\,u_{n+1}(p_{ni})=p_{ni}$.) Note that the identification of
$K_1(A)$ with $K_0(SA)$ is done in such a way that $[u_n]$
corresponds to
$$
[t\mapsto \Ad\,v_{n,t}(1\oplus 0)]-[(1\oplus0)]
$$
(\cite{Bl} 8.2.2). Since
$$
[t\mapsto \Ad\,v_{n,t}(p_{ni})]=[t\mapsto \Ad\,v_{n,t}(1\oplus0)]
-[t\mapsto \Ad\,v_{n,t}(1-p_{ni})],
$$
($*$) equals
\BE
&&[t\mapsto \Ad\,v_{n,t}(1\oplus0)]-[t\mapsto \Ad(v_{n,t}(1-p_{ni})
+v_{n+1,t}p_{ni})(1\oplus0)]\\
&&=[u_n]-[u_n(1-p_{ni})+u_{n+1}p_{ni}]=[u_{n+1}^*u_np_{ni}],
\EE
where we have used the fact that
$$
t\mapsto v_{n,t}((1-p_{ni})\oplus(1-\al(p_{ni}))+
v_{n+1,t}(p_{ni}\oplus\al(p_{ni}))
$$
is a path of unitaries from $1\oplus1$ to
$$
(u_n(1-p_{ni})+u_{n+1}p_{ni})\oplus(u_n^*(1-\al(p_{ni}))+u_{n+1}^*\al(p_{ni}))
.
$$
\end{pf}

\begin{lem}\label{3.9}
Suppose that $\tilde{\eta}(\al)=0$. Then there is a unitary $u_n\in A$
for each $n$ such that
\BE
&& \al|B_n=\Ad\,u_n|B_n,\\
&&\|\al(z_m)-\Ad\,u_n(z_m)\|<2^{-n},\ \ m\leq n,\\
&& B(u_{n+1}^*u_n,z_{ni})=0,\\
&&[u_{n+1}^*u_np_{ni}]=0,\\
&&\hat{h}_{ni}=0
\EE
for $i=1,\ldots,k_n$, where
$$
h_{ni}=\frac{1}{2\pi i}\log\,\al(z_{ni})\Ad\,u_n(z_{ni}^*).
$$
\end{lem}
\begin{pf}
Comparing with Lemma~\ref{3.7}, the newly appeared conditions are only
$$
[u_{n+1}^*u_np_{ni}]=0.
$$
We will find a unitary $w_n\in A\cap B_n'$ such that
$[w_n,z_n]=0$ and the above conditions are satisfied by replacing all $u_n$
by $u_nw_n$. With the condition $[w_{n+1},z_{n+1}]=0$, it follows that
$[w_{n+1},z_n]=0$ and that the other conditions are preserved.

From the assumption that
$$
0\lra K_1(A)\lra K_0(\Mal)\lra K_0(A)\lra 0
$$
is trivial, we have  a homomorphism 
$h_n^1:K_0(A_n)\ra K_1(A)$ for each $n$ such that
$$
\psi_n^1=h_n^1-h_{n+1}^1\chi_n^0.
$$
We only have to find a unitary $w_n\in A\cap B_n'$ such that
$[w_n,z_n]=0$ and
$$
[w_np_{ni}]=-[n,i]h_n^1(e_i),\ \ i=1,\ldots,k_n.
$$
Since $z_np_{ni}$ in $p_{ni}A_mp_{ni}\cap B_{ni}'$ for $m>n$ 
is a direct sum of elements
of the form
$$
\left(\begin{array}{cccc}
0&&&z_{n+1}^Lp_{ni}\\
1&\cdot&&\\
&\ddots&\ddots&\\
&&1&0
\end{array}\right)
$$
with $L=\pm1$, this follows immediately.
\end{pf}

\noindent
{\em Proof of (1)$\Rightarrow$(2) of Theorem 3.1.}

Under the assumption (1) we have found a sequence $\{u_n\}$ of unitaries
as in the previous lemma. Now we apply the homotopy lemma to the pair
$u_{n+1}^*u_np_{ni},z_np_{ni}$ of unitaries in $p_{ni}Ap_{ni}\cap B_{ni}'$
(\cite{BEEK} 8.1): From the conditions
\BE
&& B(u_{n+1}^*u_n,z_{ni})=0,\\
&& [u_{n+1}^*u_np_{ni}]=0
\EE
calculated in $p_{ni}Ap_{ni}\cap B_{ni}'$, which follow since
$K_*(p_{ni}Ap_{ni}\cap B_{ni}')\ra K_*(p_{ni}Ap_{ni})\ra K_*(A)$
are injective, and the condition $\|[u_{n+1}^*u_n,z_n]\|\ra0$ as $n\ra\infty$,
we obtain a continuous path $v_{ni,t}$ of unitaries in
$p_{ni}Ap_{ni}\cap B_{ni}'$ such that
$$
v_{ni,0}=p_{ni},\ \ v_{ni,1}=u_n^*u_{n+1}p_{ni}
$$
and
$$
\max_t \|[v_{ni,t},z_{ni}]\|\lra0\ {\rm as}\ n\ra\infty.
$$
Let $v_{nt}=\sum_i v_{ni,t}$, and define a continuous path $v_t$ of
unitaries for $t\in[1,\infty)$ by
\BE
v_1&=&u_1,\\
v_{n+t}&=&u_nv_{n,t},\ \ 0\leq t \leq 1
\EE
for $n=1,2,\ldots.$ Then since
$$
\max_t \|[v_{n,t},z_m]\|\lra0 \ {\rm as}\ n\ra\infty,
$$
we obtain that for any $m$,
$$
\lim_{t\rightarrow\infty}\Ad\,v_t(z_m)=\al(z_m).
$$
We also have that for $t\geq m$ and $a\in B_m$
$$
\Ad\,v_t(a)=\al(a).
$$
Thus it follows that for any $x\in A$
$$
\lim_{t\rightarrow\infty}\Ad\,v_t(x)=\al(x).
$$
This completes the proof.

\section{Main Theorem}
\setcounter{theo}{0}
\begin{prop}
If $\varphi\in\Hom(K_1(A),\Aff)$, there exists an automorphism
$\al\in\Ap$ such that $\eta(\al)$ is trivial 
and the rotation map $R_{\alpha}:K_1(\Mal)\ra \Aff$ is given by
$$
R_{\alpha}(a,b)=D(a)+\varphi(b)
$$
for some identification of $K_1(\Mal)$ with
$K_0(A)\oplus K_1(A)$.
\end{prop}

To prove this we first prepare:

\begin{lem}\label{4.2}
If $\varphi\in\Hom(K_1(A),\Aff)$, there exists an inductive system
$$
\Z^{k_1}\stackrel{\chi_1^i}{\lra}\Z^{k_2}\stackrel{\chi_2^i}{\lra}\Z^{k_3}\lra
\cdots
$$
whose limit is isomorphic to $K_i(A)$ for $i=0,1$ and homomorphisms
$h_n:\Z^{k_n}\ra\Z^{k_{n+1}}$ such that
\BE
&& |\varphi\circ\chi_{\infty,n-1}^1(e_j^{n-1})-D\circ\chi_{\infty,n}^0\circ
h_{n-1}(e_j^{n-1})|
<2^{-n+1}\el_{n-1}^{-1}D\circ\chi_{\infty,n-1}^0(e_j^{n-1}),\\
&& |h_{n-1}\circ\chi_{n-2}^1(e_j^{n-2})-\chi_{n-1}^0\circ h_{n-2}(e_j^{n-2})|
<2^{-n+3}\el_{n-2}^{-1}\chi_{n,n-2}^0(e_j^{n-2}),\\
&& |\chi_n^0(i,j)|\geq 2^{n+1}\max(|\chi_n^1(i,j)|,1),
\EE
where that $|x|<y$ for $x,y\in\Z^{k_n}$ means that $|x_i|<y_i$ for all $i$,
$(e_j^n)_j$ is the canonical basis for $\Z^{k_n}$,
$$
\el_n=\max\{[n,j]\ |\ j=1,\ldots,k_n\},
$$
and $([n,j])_j\in\Z^{k_n}$ corresponds to $[1]\in K_0(A)$.
\end{lem}
\begin{pf}
Suppose that we are given inductive systems
$$
\Z^{k_1}\stackrel{\chi_1^i}{\lra}\Z^{k_2}\lra\cdots
$$
such that the limit is isomorphic to $K_i(A)$ for $i=0,1$, and
$\chi_n^0(i,j)\geq 2^{n+1}\max(|\chi_n^1(i,j)|,1)$. By passing to a subsequence
we construct the homomorphisms $h_n$ with the required properties.

Suppose that we have constructed $h_1,\ldots,h_{n-1}$ and fixed
$\Z^{k_1},\ldots,\Z^{k_n}$. Then we compute $\el_n$ and find 
$\xi: \Z^{k_n}\ra K_0(A)$ such that
$$
|\varphi\xi_{\infty,n}^1(e_j^n)-D\xi(e_j^n)|
<2^{-n}\el_n^{-1}D\chi_{\infty,n}^0(e_j^n).
$$
This is obviously possible by the density of $\Range\, D$ and
$$
\inf_{\tau\in T_A}D\chi_{\infty,n}^0(e_j^n)(\tau)>0.
$$
Then we find an $m>n$ such that 
$\Range\, \xi\subset \Range\,\chi_{\infty,m}^0$,
and $\eta:\Z^{k_n}\ra \Z^m$ such that
$$
\begin{array}{ccccc}
\Z^{k_n}&&\stackrel{\eta}{\lra}&&\Z^m\\
&\xi\searrow&&\swarrow\chi_{\infty,m}^0&\\
&&K_0(A)&&
\end{array}
$$
is commutative. Note
\BE
&& |D\chi_{\infty m}^0\eta\chi_{n-1}^1(e_j^{n-1})-D\chi_{\infty n}^0h_{n-1}
     (e_j^{n-1})|\\
&&\leq |D\chi_{\infty m}\eta\chi_{n-1}^1(e_j^{n-1})
     -\varphi\chi_{\infty\, n-1}^1(e_j^{n-1})|
      +|\varphi\chi_{\infty\,n-1}^1(e_j^{n-1})
           -D\chi_{\infty n}^0h_{n-1}(e_j^{n-1})|\\
&&<2^{-n}\el_n^{-1}\sum_{i=1}^{k_n}D\chi_{\infty n}^0(e_i^n)|\chi_{n-1}^1(i,j)
|
     +2^{-n+1}\el_{n-1}^{-1}D\chi_{\infty\,n-1}^0(e_j^{n-1})\\
&&<2^{-n}\el_n^{-1}\sum_iD\chi_{\infty n}^0(e_i^n)\chi_{n-1}^0(i,j)
    +2^{-n+1}\el_{n-1}^{-1}D\chi_{\infty\,n-1}^0(e_j^{n-1})\\
&&<(2^{-n}\el_n^{-1}+2^{-n+1}\el_{n-1}^{-1})D\chi_{\infty\,n-1}^0(e_j^{n-1})\\
&&<2^{-n+2}\el_{n-1}^{-1}D\chi_{\infty\,n-1}^0(e_j^{n-1}).
\EE
Thus by choosing a sufficiently large $\el>m$ it follows that
$$
|\chi_{\el m}^0\eta\chi_{n-1}^1(e_j^{n-1})-
  \chi_{\el n}^0h_{n-1}(e_j^{n-1})|
<2^{-n+2}\el_{n-1}^{-1}\chi_{\el,n-1}^0(e_j^{n-1}).
$$
By taking $\Z^{k_{\el}}$ for $\Z^{k_{n+1}}$ and $\chi_{\el m}^0\eta$
for $h_n$, this completes the proof.
\end{pf}

\noindent
{\em Proof of Proposition 4.1}

By the previous lemma we have the following diagram:
$$
\begin{array}{cccccccccc}
\lra&\Z^{k_n}&\stackrel{\chi_n^1}{\lra}&\Z^{k_{n+1}}&\lra&&\cdots&&\lra
      &K_1(A)\\
&&\searrow h_n&&\searrow h_{n+1}&&&&&\\
\lra&\Z^{k_n}&\stackrel{\chi_n^0}{\lra}&\Z^{k_{n+1}}&\lra&\Z^{k_{n+2}}&\lra
  &\cdots&\lra&K_0(A)
\end{array}
$$
with the specified properties. Accordingly we construct an increasing
sequence $\{A_n\}$ of T algebras such that
\BE
A_n&=&B_n\otimes C(\T),\\
B_n&=&\oplus_{i=1}^{k_n}B_{ni},\\
B_{ni}&\cong& M_{[n,i]}
\EE
and the embeddings of $A_n$ into $A_{n+1}$ are in the standard form. 
By Elliott's
theory \cite{El}, we identify $\ol{\cup_{n=1}^{\infty}A_n}$ with $A$.

Define $\psi_n^0:K_1(A_n)\ra K_0(A_{n+2})$ by
$$
\psi_n^0=h_{n+1}\chi_n^1-\chi_{n+1}^0h_n.
$$
By the properties specified in Lemma~\ref{4.2} we have that
$$
|\psi_n^0(i,j)|<2^{-n+1}\el_n^{-1}\chi_{n+2,n}^0(i,j).
$$
Then by Lemma~\ref{3.4} (and its proof) we find a unitary
$w_{nj}\in B_{n+2}\cap B_n'$ such that
\BE
&& w_{nj}=w_{nj}p_{nj}+1-p_{nj},\\
&& \|\Ad\,w_{nj}(z_{nj})-z_{nj}\|\leq 3\pi 2^{-n+1},\\
&& B(w_{nj},z_{nj})=-[n,j]\psi_n^0(e_j^n).
\EE
(Because $z_{nj}$ in $B_{n+2,i}\otimes C(\T)$ is a direct sum of elements
of the form as in the proof of Lemma~\ref{3.4} such that the matrix sizes
$M_s$ are at least $2^{2n}$; hence the error introduced by choosing $N_s$
in that proof will be of the order $2^{-2n}$.)
If $w_n$ denotes $w_{n1}w_{n2}\cdots w_{nk_n}$, then we have that
\BE
&& w_n\in B_{n+2}\cap B_n',\\
&& \|\Ad\,w_n(z_n)-z_n\|\leq 3\pi 2^{-n+1},\\
&& B(w_n,z_{nj})=-[n,j]\psi_n^0(e_j^n),\\
&& [w_np_{nj}]=0.
\EE
We define the following two automorphisms $\beta_0,\beta_1$ of $A$ by
\BE
\beta_0&=&\lim_{n\rightarrow\infty}\Ad(w_2w_4\cdots w_{2n}),\\
\beta_1&=&\lim_{n\rightarrow\infty}\Ad(w_1w_3\cdots w_{2n-1}).
\EE
To show the limits exist, note that $[w_m,w_n]=0$ if $|m-n|\geq2$ and
the limits obviously exist on $\cup_{n=1}^{\infty}B_n$. Since 
$\Ad(w_nw_{n+2}\cdots w_{n+2k})(z_n)$ in $A_{n+2k+2}$ is a direct sum of
elements of the form
$$
\left(\begin{array}{cccc}
0&&&z_{n+2k+2}^L\\  1&\cdot&&\\  &\ddots&\ddots&\\ &&1&0
\end{array}\right)
$$
with $L=\pm1$, we have that
$$
\|\Ad(w_n\cdots w_{n+2k}w_{n+2k+2})(z_n)-\Ad(w_n\cdots w_{n+2k})(z_n)\|
<3\pi 2^{-(n+2k+1)}.
$$
Then it also follows that the limits exist on $z_1,z_2,\ldots$. Since the
same reasoning applies to the inverses, we have shown that $\beta_0,\beta_1$
exist as automorphisms.

Now we shall show that the product $\beta_0\beta_1$ satisfies the required
properties.

By \cite{KK1}, 2.4 the extension $\eta_1(\beta_i)$
$$
0\lra K_1(A)\lra K_0(M_{\beta_i})\lra K_0(A)\lra 0
$$
is trivial for $i=0,1$ and the extension $\eta_0(\beta_i)$
$$
0\lra K_0(A)\lra K_1(M_{\beta_i})\lra K_1(A)\lra0
$$
is given as the inductive limit of
$$
\begin{array}{ccccccccccc}
0&\lra &\Z^{k_n}&\lra&\Z^{k_n}&\oplus&\Z^{k_n}&\lra&\Z^{k_n}&\lra&0\\
&&{\scriptstyle \chi_{n+2,n}^0}\downarrow&&\downarrow&
  {\scriptstyle \psi_n^0}\swarrow&\downarrow&&
  \downarrow&&\\ 
0&\lra&\Z^{k_{n+2}}&\lra&\Z^{k_{n+2}}&\oplus&\Z^{k_{n+2}}&\lra&\Z^{k_{n+2}}&
    \lra&0\\
&&\downarrow&&\downarrow&\swarrow&\downarrow&&\downarrow&&
\end{array}
$$
with $n\equiv i$ (mod 2). Hence $\eta_1(\beta_0\beta_1)=\eta_1(\beta_0)
+\eta_1(\beta_1)=0$. We will compute $\eta_0(\beta_0)+\eta_0(\beta_1)$
below.

Define 
$$
E=\{(x,y)\in K_1(M_{\beta_0})\oplus K_1(M_{\beta_1})\,|\,q(x)=q(y)\}
/\{(a,-a)\,|\,a\in K_0(A)\}.
$$
If $g\in K_1(A)$ is the image of $x_{2n+1}\in\Z^{k_{2n+1}}$, define
$\eta_n:\Range\,\chi_{\infty\, 2n+1}^1 
\ra K_1(M_{\beta_0})\oplus K_1(M_{\beta_1})$ by
$$
\eta_n(g)=(h_{2n+1}(x_{2n+1}),x_{2n+2})\oplus(0,x_{2n+1}),
$$
where the right hand side should be regarded as an element of
$K_1(M_{\beta_0})\oplus K_1(M_{\beta_1})$. Then
\BE
&&\eta_{n+1}(g)-\eta_n(g)\\
&&=(h_{2n+3}(x_{2n+3})-\psi_{2n+2}^0(x_{2n+2})
    -\chi_{2n+4,2n+2}^0h_{2n+1}(x_{2n+1}),\,0)\\
&&\ \ \oplus
    (-\psi_{2n+1}^0(x_{2n+1}),0)\\
&&=(\chi_{2n+3}^0h_{2n+2}(x_{2n+2})-\chi_{2n+4,2n+2}^0
     h_{2n+1}(x_{2n+1}),\,0)\\ 
&&\ \  \oplus(-h_{2n+2}(x_{2n+2})+\chi_{2n+2}^0h_{2n+1}(x_{2n+1}),0).
\EE
Thus $(\eta_n)$ gives a well-defined homomorphism $\eta:K_1(A)\ra E$
such that $q\eta=\id$. This shows that $\eta_0(\beta_0\beta_1)=0$.

Let $u_n=w_nw_{n-2}\cdots$.
We take a path $v(t)$ of unitaries in $A\otimes M_2$ from $z_{nj}$ to
$\beta_0(z_{nj})$ by composing the following two paths for even $m\geq n$:
$$
v_1(t)=R_t(1\oplus u_m)R_t^{-1}(z_{nj}\oplus1)R_t(1\oplus u_m^*)R_t^{-1}
$$
and a short path $v_2$ from $\Ad\,u_m(z_{nj})$ to $\beta_0(z_{nj})$.
For $\tau\in T_A$ we want to compute
$$
\frac{1}{2\pi i}\int_0^1\tau(\dot{v}(t)v(t)^*)dt.
$$
We know the contribution from $v_1$ is zero and the contribution from $v_2$
is given by
\BE
&&\lim_{k\rightarrow\infty}\tau(B(w_{m+2}^*w_{m+4}^*\cdots w_{m+2k}^*,z_{nj}))
\\
&&=\lim\tau(\sum_{i=1}^k\chi_{\infty\,2m+2i+2}^0\psi_{m+2i}^0\chi_{m+2i,n}^1
    (e_j^n)).
\EE
Thus we obtain that
$$
R_{\beta_0}([v])=\sum_{i=1}^{\infty}D\chi_{\infty\,2m+2i+2}^0\psi_{m+2i}^0
    \chi_{m+2i,n}^1(e_j^n).
$$
A similar computation applies to $\beta_1$.
For an odd $n$ we let $m=n-1$ for computing $r_0=R_{\beta_0}([v])$ and
let $m=n$ for computing the corresponding $r_1$, and obtain that
\BE
r_0+r_1&=&\sum_{i=1}^{\infty}D\chi_{\infty\,n+i+2}^0\psi_{n+i}^0\chi_{n+i,n}^1
            (e_j^n)\\
&=&\sum_{i=1}^{\infty}\left(D\chi_{\infty\,n+i+2}^0h_{n+i+1}\chi_{n+i+1,n}^1
       (e_j^n)-D\chi_{\infty\,n+i+1}^0h_{n+i}\chi_{n+i,n}^1(e_j^n)\right)\\
&=&\varphi\chi_{\infty\,n}^1(e_j^n)-D\chi_{\infty\,n+1}^0h_{n+1}\chi_n^1
   (e_j^n).
\EE
Under the identification of $K_1(M_{\beta_0\beta_1})$ with
$K_0(A)\oplus K_1(A)$ specified above, the above element corresponds
to $(-h_{n+1}\chi_n^1(e_j^n),[z_{nj}])$. This implies that
$R_{\beta_0\beta_1}$ satisfies the required properties.

\medskip
Let $Q$ be the homomorphism of $\OExt(K_1(A),K_0(A))$ into
$\Ext(K_1(A),K_0(A))$ defined by $[(E,R)]\mapsto [E]$. Then $\ker Q$
is the subgroup of the isomorphism classes of $(E_0,R_{\varphi})$ where
$E_0$ is the trivial extension $K_1(A)\oplus K_0(A)$, and
$R_{\varphi}:E_0\ra\Aff$ is determined by $\varphi\in\Hom(K_1(A),\Aff)$
as in the previous proposition:
$$
R_{\varphi}:(a,b)\mapsto D(a)+\varphi(b).
$$

\begin{prop}\label{4.3}
The following sequences of abelian groups are exact:
\BE
&& 0\lra\ker Q\lra\OExt(K_1(A),K_0(A))\stackrel{Q}{\lra}\Ext(K_1(A),K_0(A))
   \lra0,\\
&&0\lra\Hom(K_1(A),\ker D)\lra\Hom(K_1(A),K_0(A))\\
&&\ \ \ \ \ \lra\Hom(K_1(A),\Aff)
  \lra\ker Q\lra0.
\EE
\end{prop}
\begin{pf}
For the first sequence we only have to show that $Q$ is surjective.
Given an extension
$$
0\lra K_0(A)\lra E\lra K_1(A)\lra 0,
$$
we regard $K_0(A)$ as a subgroup of $E$ and have to extend $D:K_0(A)\ra\Aff$
to a homomorphism $R:E\ra\Aff$. This can be done step by step by using
the fact that \Aff\ is divisible.

For the second sequence we only have to show that $(E_0,R_{\varphi})$
and $(E_0,R_{\psi})$ are isomorphic if and only if 
$\varphi=\psi+D\circ h$ for some $h\in\Hom(K_1(A),K_0(A))$. 
This follows because an
isomorphism $\mu:E_0\ra E_0$ is given by
$$
\mu:(a,b)\mapsto (a+h(b),b)
$$
for some $h\in\Hom(K_1(A),K_0(A))$ with $R_{\psi}\circ \mu=R_{\varphi}$.
\end{pf}

\begin{theo}\label{4.4}
Let $A$ be a simple unital AT algebra of real rank zero, \Ap\ the group of
approximately inner automorphisms of $A$, and \AI\ the group of 
asymptotically inner automorphisms of $A$. Then \AI\ is a normal subgroup
of \Ap\ and the quotient $\Ap/\AI$ is isomorphic to
$$
\OExt(K_1(A),K_0(A))\oplus\Ext(K_0(A),K_1(A))
$$
with isomorphism induced by $\tilde{\eta}$.
\end{theo}
\begin{pf}
Before Theorem~\ref{3.1} we have described the homomorphism
$$
\tilde{\eta}:\Ap\ra\OExt(K_1(A),K_0(A))\oplus\Ext(K_0(A),K_1(A)),
$$
and showed in \ref{3.1} that $\ker \tilde{\eta}=\AI$. By 3.1 of \cite{KK1}
we have shown that $\eta=(\eta_0,\eta_1)=(Q\tilde{\eta}_0,\eta_1)$
is surjective onto $\Ext(K_1(A),K_0(A))\oplus\Ext(K_0(A),K_1(A))$.
By Proposition 4.1 we know that $\Range\, \tilde{\eta}$ contains
$\ker Q$, which shows that $\tilde{\eta}$ is surjective.
This completes the proof.
\end{pf}

\begin{ex}\label{4.5}
If $A$ is the irrational rotation \cstar\ generated by unitaries
$u,v$ with $uvu^*v^*=e^{2\pi i\theta}1$ for some irrational
number $\theta\in(0,1)$, then $A$ is a simple unital AT algebra
of real rank zero by \cite{EE}, and
$K_i(A)\cong \Z^2$ and hence $\Ext(K_i(A),K_{i+1}(A))=0$. But
since $A$ has only one tracial state and $\Range\,D=\Z+\theta\Z$, 
it follows that $\Hom(K_1(A),\Aff)\cong \R^2$ and $\OExt(K_1(A),K_0(A))
\cong \R^2/(\Z+\theta\Z)^2$ which is isomorphic to $\Ap/\AI$.
Note also that \HI=\Ap\ in this case since the natural $\T^2$ action on
$A$ exhausts all \OExt.
\end{ex}

\medskip
\small
\begin{flushright}
Department of Mathematics, Hokkaido University, Sapporo 060 Japan\\
Department of Mathematics, University of Nevada, Reno NV 89557 USA
\end{flushright}
\end{document}